\documentclass{article}
\usepackage{amsmath}
\usepackage{amssymb}
\usepackage[utf8]{inputenc}
\title{ERRATUM TO ``ALMOST MENGER PROPERTY IN BITOPOLOGICAL SPACES"}
\author{Santanu Acharjee$^{1,\dagger}$, Kabindra Goswami$^2$\\
$^{1}$Department of Mathematics\\
 Gauhati University\\
 Guwahati-781014, Assam, India\\
$^{2}$Department of Mathematics\\
Goalpara College\\
Goalpara-783101, Assam, India\\
e-mails: $^{1}$sacharjee326@gmail.com,$^{2}$kabindragoswami@gmail.com\\ 
$^\dagger$Orcid id: 0000-0003-4932-3305\\
$^\dagger$Corresponding author: Santanu Acharjee}

\date{}

\begin{document}
\maketitle
{\bf Abstract:} In this erratum, theorem 3.2 and theorem 4.1  of \"{O}z\c{c}a\u{g} and Eysen [S. \"{O}z\c{c}a\u{g}  and A.E. Eysen,  Almost Menger property in bitopological spaces, Ukrainian Math. J., {\bf 68}, No 6, 950-958 (2016)] are proven to be incorrect. An example is provided to disprove them and thus, correct versions of  the theorems are restated.\\

{\bf UDC} 517.9

\section{Introduction} Bitopological space  \cite{Kelly} is one the areas of topology which has attracted attentions of researchers. It was introduced by Kelly \cite{Kelly} in 1963. Later, many researchers have worked in this area. Recently,  \"{O}z\c{c}a\u{g} and Eysen \cite{Ozcag} studied  almost Menger property in bitopological spaces and proved several results. But, theorem 3.2 and theorem 4.1 of \cite{Ozcag} are  found to be incorrect. Thus,  suitable example is provided to disprove their claims and then, the correct versions of the theorem 3.2 and theorem 4.1 are reestablished.\\

\section{Corrigendum}

In this section,  theorem 3.2 and theorem 4.1 of \cite{Ozcag} are proven to be incorrect. Later,  correct versions of the  theorems are restated.\\

\textbf{Definition 1.} (\cite {Ozcag},  Definition 3.2.) 
Let $(X, \tau_1, \tau_2)$ be a bitopological space. A set $A\subseteq X$ is called $(i, j)$-regular
open ($(i, j)$-regular closed)($i \neq j, i, j = 1, 2$) if $ A = Int_{\tau_i} Cl_{\tau_j} (A )$  $(A=Cl_{\tau_i} Int_{\tau_j} (A))$. The set $A$ is called pairwise regular open (pairwise regular closed) if it is both $(i, j)$-regular open and $(j, i)$-regular open
($(i, j)$-regular closed and $(j, i)$-regular closed).\\

\textbf{Theorem 1.} (\cite {Ozcag},  Theorem 3.2.) 
A bitopological space $(X, \tau_1, \tau_2)$ is $(i, j)$-almost Menger if and only if, for each sequence $\langle \mathcal{U}_n :n \in \mathbb{N} \rangle$ of covers of $X$ by $(i, j)$-regular open sets, there exists a sequence $\langle \mathcal{V}_n :n \in \mathbb{N} \rangle$ of finite families such that, for any $n \in \mathbb{N}$,
\[\mathcal{V}_n \subseteq \mathcal{U}_n \; \text{and} \; X=\bigcup_{n \in \mathbb{N}}(\bigcup_{V \in \mathcal{V}_n}  Cl_{\tau_j}(V))\]\\

\textbf{Theorem 2.} (\cite {Ozcag},  Theorem 4.1.) 
$(X, \tau_1, \tau_2)$ is an $(i, j)$-almost $\gamma$-set if and only if, for any sequence $\langle \mathcal{U}_n :n \in \mathbb{N} \rangle$ of
$\tau_i$-$\omega$-covers of $X$ by $(i, j)$-regular open sets, there exists a sequence $\langle V_n :n \in \mathbb{N} \rangle$ such that, for all $n \in \mathbb{N}$, $V_n \in \mathcal{U}_n$ and the set $\{V_n :n \in \mathbb{N}\}$ is an $(i, j)$-almost $\gamma$-cover for $X$.\\

Now,  the following example is considered to disprove above two theorems.\\ 

\textbf{Example 1.} Let us consider a bitopological space $(X,\tau_1,\tau_2)$, where $X=\{a,b,c, d\}$, $\tau_1=\{\emptyset, \{a\}, \{b\}, \{a,b\}, X\}$ and $\tau_2=\{\emptyset, \{a\}, \{c\}, \{a,c\}, X\}$. At first, we consider $U_1=\{a,b\}$, then clearly $U_1\in \tau_1$. Thus,  $Int_1(Cl_2(U_1)) = \{a,b\}=U_1$. Again, if we consider $V_1=\{a\}$, then clearly $V_1\in \tau_1$ and  $Int_1(Cl_2(V_1)) = \{a,b\}\neq V_1$. So, $U_1$ is $(1,2)$-regular open but $V_1$ is not $(1,2)$-regular open.
Thus, it is not true in general that if $U\in \tau_1$, then $U$ is $(1,2)$-regular open in $(X,\tau_1,\tau_2)$.\\

Again, we consider $U_2= \{a,c\}$. Then, $U_2\in \tau_2$ and $Int_2(Cl_1(U_2)) = \{a,c\}=U_2$. Again, if we consider $V_2= \{a\}$. Then, $V_2\in \tau_2$ and $Int_2(Cl_1(V_2)) =\{a,c\} \neq V_2$. So, $U_2$ is $(2,1)$-regular open but $V_2$ is not $(2,1)$-regular open.
Thus, it is not true in general that if $U\in \tau_2$, then $U$ is $(2,1)$-regular open in $(X,\tau_1,\tau_2)$.\\

To prove sufficient part of theorem 3.2, \"{O}z\c{c}a\u{g} and Eysen \cite {Ozcag} considered $\langle \mathcal{U}_n : n\in \mathbb{N} \rangle$. Later for any n, \"{O}z\c{c}a\u{g} and Eysen \cite {Ozcag} considered $ \mathcal{U}_n^{'} = \{Int_1(Cl_2(U))  : U\in \mathcal{U}_n \}$. According to them, $\langle \mathcal{U}_n^{'} : n\in \mathbb{N} \rangle$ is a sequence of covers of $X$ by $(1,2)$-regular open sets. But,  example 1 shows that the claim made in the proof of the sufficient part of  \cite {Ozcag} is not true in general. It is not true in general that if $U\in \tau_1$, then $U$ is $(1,2)$-regular open. Thus, the proof of the sufficient part is incorrect. Hence,   theorem 3.2 of \cite{Ozcag} is rectified as given below.\\

\textbf{Theorem 3.} 
If a  bitopological space $(X, \tau_1, \tau_2)$ is $(i, j)$-almost Menger, then  for each sequence $\langle \mathcal{U}_n :n \in \mathbb{N} \rangle$ of covers of $X$ by $(i, j)$-regular open sets, there exists a sequence $\langle \mathcal{V}_n :n \in \mathbb{N} \rangle$ of finite families such that, for any $n \in \mathbb{N}$,
\[\mathcal{V}_n \subseteq \mathcal{U}_n \; \text{and} \; X=\bigcup_{n \in \mathbb{N}}(\bigcup_{V \in \mathcal{V}_n}  Cl_{\tau_j}(V))\]\\

Using the similar claim made in example 1,  theorem 4.1 of \cite {Ozcag} is disproved. Thus, the correct version of  theorem 4.1 of \cite{Ozcag} can be stated as below.\\ 

\textbf{Theorem 4.} 
If $(X, \tau_1, \tau_2)$ is an $(i, j)$-almost $\gamma$-set, then  for any sequence $\langle \mathcal{U}_n :n \in \mathbb{N} \rangle$ of $\tau_i$-$\omega$-covers of $X$ by $(i, j)$-regular open sets, there exists a sequence $\langle V_n :n \in \mathbb{N} \rangle$ such that, for all $n \in \mathbb{N}$, $V_n \in \mathcal{U}_n$ and the set $\{V_n :n \in \mathbb{N}\}$ is an $(i, j)$-almost $\gamma$-cover for $X$\\

\section{Conclusion}
In this paper,  theorem 3.2 and theorem 4.1 of \cite {Ozcag} were proven to be incorrect. An example was provided to disprove the theorems and later, the correct versions of the theorems were restated. It can be assumed that this paper may be useful in  research of bitopological selection principles in future.\\

{\bf Declarations}

{\bf Funding: } No funds, grants, or other support was received.

{\bf Competing interests} The authors declare that there is no competing interest.

\end{document}